\documentclass[12pt]{article}
\usepackage[pctex32]{graphics}
\usepackage{amsthm,amsmath,amssymb,amscd,verbatim,epsfig}
\usepackage{amssymb}
\usepackage{graphicx}
\newtheorem{Theorem}{Theorem}
\newtheorem{Lemma}{Lemma}

\newcommand{\beq}{\begin{equation}}
\newcommand{\eeq}{\end{equation}}

\date{}
\newcommand{\al}{\alpha}

\newcommand{\f}{\frac}

\newcommand{\va}{\varepsilon}

\newcommand{\sq}{$\square$}

\begin{document}

\title{The\ lives\ of\ period-3\ orbits\ for\ some\ quadratic\ polynomials}
\author{Bau-Sen Du \\ [.3cm]
Institute of Mathematics \\
Academia Sinica \\
Taipei 11529, Taiwan \\
dubs@math.sinica.edu.tw \\}
\maketitle


\section{Introduction}
In the past 25 years, there have been a number of papers {\bf {\cite{Bech, Gor, Li1, Mir, Saha}}} dealing with the bifurcation of period-3 orbits of the logistic map $F_\mu(x) = \mu x(1-x)$.  Although it is well known that, for quadratic polynomial such as the logistic map $F_\mu$, as the parameter $\mu$ increases from 0, new periodic orbits are created through tangent bifurcations or period-doubling bifurcations and never die once they have been born, it is less known that such phenomena do not hold for general quadratic polynomials.  Given the pedagogical value of the logistic map which has been used by many texts {\bf {\cite{Dev, Ela, Rob}}} to illustrate various aspects of nonlinear dynamics such as equilibrium, periodicity, stability, structural stability, bifurcations, chaos, etc., we would like to present some examples of quadratic polynomials which show different phenomena.  Furthermore, when one deals with the bifurcation of periodic orbits, one often needs to use the implicit function theorem.  Since, as noted in {\bf {\cite{Dev}}}, most beginning students do not appreciate the power of the implicit function theorem when they encounter it in their first analysis course, it may be a good idea to give a nontrivial application of the theorem to a concrete example to help demonstrate its usefulness.  The purpose of this note is, on the one hand, to give a full description of the lives of period-3 orbits of the family $f_\al(x) = 1 - \al x^2$ via a nontrivial application of the implicit function theorem and, on the other, to use it to show that some one-parameter families of quadratic polynomials (such as $T_{a,c}(x) = a - c(1+x^2)$ with suitably chosen $a > 0$) possess such peculiar pathologies as "bubbles" (i.e., periodic orbits live a finite life) and/or point bifurcations (i.e., periodic orbits die right on birth) {\bf {\cite{Du3, Li2, Nus}}}.  

\section{The immortal lives of period-3 orbits for $S_{a,c}(x) = a - cx^2$}
The following result whose proof uses only straightforward calculations is taken from {\bf {\cite{Du1}}}.  

\begin{Lemma}
\label{lemma1}
For $\al > 0$, let $f_\al(x) = 1 - \al x^2$.  Then
\begin{itemize}
\item[(a)]
when $\al < 7/4$, $f_\al(x)$ has no period-3 points,

\item[(b)]
when $\al = 7/4$, $f_\al(x)$ has 3 period-3 points (which are zeros of $343x^3-98x^2-252x+8$),

\item[(c)]
when $\al > 7/4$, $f_\al(x)$ has 6 period-3 points.
\end{itemize}
\end{Lemma}

{\bf Proof.}
Let $p_\al(x) = 1 - x - \al x^2$.  Then $1-\al x^2 = p_\al(x)+x$ and  $f_\al^3(x) - x = 1-\al[1-\al(1-\al x^2)^2]^2 -x = 1-\al[1-\al(p_\al(x) + x)^2]^2 -x$ $= p_\al(x)h(\al, x)$, where $h(\al, x)$ $=1-\al^3p_\al^3(x) + (-4\al^3 x+2\al^2)p_\al^2(x) +(-4\al^3x^2 +6\al^2 x-\al)p_\al(x)+4\al^2 x^2-2\al x$
$=\al^6x^6 - \al^5x^5 + (-3\al^5+\al^4)x^4+(2\al^4-\al^3)x^3+(3\al^4-3\al^3+\al^2)x^2
+(-\al^3+2\al^2-\al)x$
$-\al^3+2\al^2-\al+1$.

We note that $h(\al_0, x_0) = 0$ if and only if $x_0$ is a period-3 point of $f_{\al_0}(x)$.  Thus, if $h(\al_0, x)$, as a polynomial in $x$,  has a real zero $x_0$, then $x_0, f_{\al_0}(x_0)$, and $f_{\al_0}^2(x_0)$ are three distinct real zeros of $h(\al_0, x)$ and, since any real polynomial of odd degree must have at least one real zero, all 6 zeros of $h(\al_0, x)$ are real.

Assume that $h(\al_0, x_0) = 0$ and $\f {\partial}{\partial x} h(\al_0, x_0)= 0$.   Then the multiplicity of $x_0$ is at least 2 and so, $h(\al_0, x)$ must be a complete square.  Let $h(\al_0, x) = (ax^3+bx^2+cx+d)^2$.  We may assume that $a > 0$.  By comparing the coefficients of $x^6, x^5, x^4$, and $x^3$ on both sides respectively, we obtain that $a = \al_0^3$, $b = -\al_0^2/2$, $c = - 3\al_0^2/2+ 3\al_0/8$, and $d = {\al_0}/4- 5/{16}$.  By comparing the coefficients of $x$ on both sides, we obtain that $\al_0 = 7/4$ is the only possible solution.  By checking with the coefficients of $x^2$ and constant terms on both sides, we obtain that $\al_0 = 7/4$ is really the only solution and $h(7/4, x) = (ax^3+bx^2+cx+d)^2 =(1/64)^2[343x^3-98x^2-252x+8)]^2$.  Conversely, if $\al_0 = 7/4$, then it follows from the above that $h(7/4,x) = (1/64)^2[343x^3-98x^2-252x+8)]^2$.  So, if $x_0$ is any real zero of the polynomial $343x^3-98x^2-252x+8$, then, we obtain that $h(\al_0, x_0) = 0$ and $\f {\partial}{\partial x} h(\al_0, x_0)= 0$.  This proves (b).

Now, there are 3 changes in signs of the coefficients of $h(2, x)$.  So, $h(2, x)=0$ has at least one and hence 6 real solutions.  By implicit function theorem, we can continue each of these 6 solutions further from $\al = 2$ as long as $\f {\partial}{\partial x} h(\al, x) \ne 0$  and this inequality holds as long as $\al > 7/4$.  This proves (c).  Again, by implicit function theorem, we have (a) since $h(0, x) (\equiv 1)$ has no real zeros.
\hfill\sq

{\bf{Remark.}}
Let $F_\mu(x) = \mu x(1-x)$ and let $G_\mu(x) = 1 - [(\mu^2-2\mu)/4]x^2$.  Then $F_\mu(x)$ is topologically conjugate to $G_\mu(x)$ through $h_\mu(x) = 1/2 + (\mu/4 - 1/2)x$ $:(F_\mu \circ h_\mu)(x) = \mu/4 -\mu(\mu/4 - 1/2)^2x^2 = (h_\mu \circ G_\mu)(x)$.  Since topological conjugacy preserves periodic points of same periods {\bf {\cite{Dev, Ela, Rob}}}, $F_\mu(x)$ has a period-3 bifurcation, by Lemma~\ref{lemma1}, at $\mu = 1 + 2 \sqrt{2}$ and exist for all $\mu \ge 1 + 2 \sqrt{2}$.

For $a \ne 0$, let $S_{a,c}(x) = a - cx^2$ and $R_a(x) = ax$.  Then $(S_{a,c} \circ R_a)(x) = a - a^2cx^2 = a(1 - acx^2) = (R_a \circ f_{ac})(x)$.  That is, $S_{a,c}(x)$ is topologically conjugate to $f_{ac}(x)$ through $R_a(x)$.  Therefore, by Lemma~\ref{lemma1}, we immediately obtain the following result which shows that, for the family $S_{a,c}(x) = a - cx^2$ with either $a$ or $c$ fixed, once the period-3 orbits are born they live for eternity.

\bigskip
\begin{Theorem}
\label{theorem1}
For $ac > 0$, let $S_{a,c}(x) = a - c x^2$.  Then
\begin{itemize}
\item[(a)]
when $ac < 7/4$, $S_{a,c}(x)$ has no period-3 orbits,

\item[(b)]
when $ac = 7/4$, $S_{a,c}(x)$ has exactly one period-3 orbit,

\item[(c)]
when $ac > 7/4$, $S_{a,c}(x)$ has exactly two period-3 orbits.
\end{itemize}
\end{Theorem}

In the above theorem, by replacing $a$ with $a-bc$, we obtain the new map $T_{a,b,c}(x) = a - c(b+x^2)$.  The period-3 orbits of $T_{a,b,c}(x)$ live very different lives from that of $S_{a,c}(x)$ as shown below (for simplicity, we choose $b = 1$ and let $T_{a,c}(x) = a - c(1+x^2)$).

\section{The finite lives of period-3 orbits for $T_{a,c}(x) = a - c(1+x^2)$ with fixed $a > \sqrt{7} \approx 2.64575$}
Let $g_c$ be a one-parameter family of continuous maps from the real line into itself.  We say that $g_c$ has a bubble of period-$n$ orbit on $[c_1, c_2]$ if there exist two continuous real-valued functions $x_i(c)$, $i = 1, 2$, on $[c_1, c_2]$ such that $x_1(c_i) = x_2(c_i)$, $i = 1,2$, $x_1(c) < x_2(c)$ for all $c_1 < c < c_2$ and each $x_j(c)$, $j = 1, 2$, is a period-$n$ point of $g_c$.  The following result is taken from {\bf {\cite{Du3}}} which is an easy consequence of Theorem~\ref{theorem1}.  Note that $(a-c)c = -(c-a/2)^2 + a^2/4$.  So, if $a > \sqrt{7}$, then $(a-c)c \ge 7/4$ if and only if $a/2 - \sqrt{a^2-7}/2 \le c \le a/2 + \sqrt{a^2-7}/2$.  

\begin{figure}
\centerline{\epsfig{file=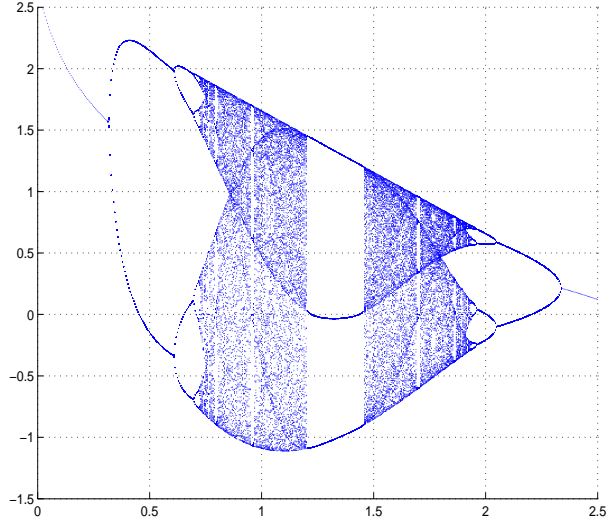,width=8cm,height=7cm}}
\caption{The bifurcation diagram with bubbles for $T_{2.658,c}(x) = 2.658 - c(1+x^2)$.}
\end{figure}

\begin{Theorem}
\label{theorem2}
For $a > 0$ and $c > 0$, let $T_{a,c}(x) = a - c(1+ x^2)$.  Then for any fixed $a > \sqrt{7}$, $T_{a,c}(x)$ has a bubble of period-3 orbit on $[a/2 - \sqrt{a^2-7}/2, a/2 + \sqrt{a^2-7}/2]$ (see Figure 1 for $a = 2.658$. Note that in the middle of the figure, the three visible curves represent the attracting period-3 orbit of $T_{2.658,c}(x)$ while the other period-3 orbit is repelling and invisible).
\end{Theorem}

{\bf{Remark.}}
Since $(a-c)c = -(c-a/2)^2 + a^2/4$, it follows from Theorem~\ref{theorem1} that, when $0 < a < \sqrt{7}$, $T_{a,c}(x)$ has no period-3 orbits.  However, it may still have bubbles of periodic orbits of some other periods.  For example, if $\sqrt{3} < a < \sqrt{7}$, then since $f_\al(x) = 1-\al x^2$ has a period doubling bifurcation from fixed point to period-2 points at $\al = 3/4$, $T_{a,c}(x)$ has a bubble of period-2 orbit on $[a/2 - \sqrt{a^2-3}/2, a/2 + \sqrt{a^2-3}/2]$ (see Figure 2 for $a = 2.35$).

\section{The tragic life of period-3 orbit for $T_{\sqrt{7},c}(x) = \sqrt{7} - c(1+x^2)$}
Let $g_c$ be a one-parameter family of continuous maps from the real line into itself.  Assume that there exist two positive numbers $\delta$ and $\va$ such that $g_c$ has a periodic point $p$ of some period $n$ for $c = c_0$, but no periodic point of same period in $(p-\va, p+\va)$ for every $c$ in $(c_0-\delta $, $c_0) \cup (c_0$, $c_0+\delta )$.  Then we say that $g_c$ has a point bifurcation of period-$n$ orbit at $c = c_0$.  The following result is taken from {\bf {\cite{Du3}}} which is an easy consequence of Theorem~\ref{theorem1}.  Note that $(a-c)c = -(c-a/2)^2 + a^2/4$.  So, if $a = \sqrt{7}$, then $(a-c)c \ge 7/4$ if and only if $c = a/2 = {\sqrt{7}}/2$.  

\begin{figure}
\centerline{\epsfig{file=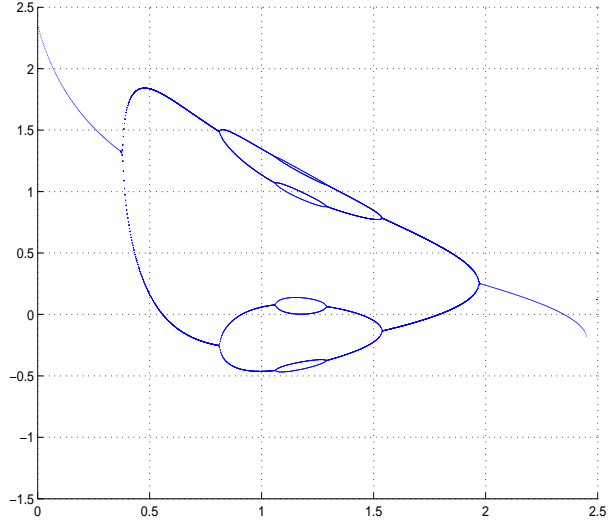,width=8cm,height=7cm}}   
\caption{The bifurcation diagram with bubbles for $T_{2.35,c}(x) = 2.35 - c(1+x^2)$.}
\end{figure}

\begin{Theorem}
$T_{\sqrt{7},c}(x)=\sqrt{7} - c(1+x^2)$ has a point bifurcation of period-3 orbit at $c = {\sqrt{7}}/2$.
\end{Theorem}

{\bf{Remark.}}
In practice, the point bifurcations of periodic orbits are very difficult to detect.  Even 
when we encounter one, we may not recognize it.  We urge the readers to experiment numerically with the family $g_c(x) = x^3 - 2x + c$ which has been shown to have a point bifurcation of period-3 orbit at $c = ?$ \, (see {\bf {\cite{Du3, Li2}}} for an answer).

\section{Some comments}
For the family of quadratic polynomials $f_\alpha(x) = 1 - \alpha x^2$, it is well-known that once a period-$n$ orbit is born it lives for eternity.  Since the family $S_{a,c}(x) = a - cx^2$ is topologically conjugate (see the paragraph before Theorem~\ref{theorem1}) to the family $f_{ac}(x) = 1 - acx^2$ when $a \ne 0$, we see that this property also holds for $S_{a,c}(x)$ for all fixed $a > 0$.  On the other hand, note that the family $T_{a,c}(x) = S_{a-c,c}(x)$ is topologically conjugate to the family $f_{(a-c)c}(x)$ and the subscript $(a-c)c$ of $f_{(a-c)c}(x)$ has a maximum (while holding $a > 0$ fixed).  This means that when we increase the value of $c$ (while holding $a > 0$ fixed) the subscript $(a-c)c$ of $f_{(a-c)c}(x)$ increases to a maximum and then decreases back to where it starts.  Consequently, if a period-$n$ orbit is born for $T_{a,c}(x)$, it either dies right on birth (point bifurcation) or lives a finite life (bubble).  So, the lives of period-$n$ orbits of $S_{a,c}(x)$ are quite different from those of $T_{a,c}(x)$.

\end{document}